\documentclass{amsart}

\usepackage{amsmath, amsfonts, amsthm, amssymb, amscd}

\usepackage{manfnt}

\usepackage[all]{xy}

\usepackage[normalem]{ulem}
\usepackage{bm,bbm}

\usepackage[bookmarks=true,
            bookmarksnumbered=true, breaklinks=true,
            pdfstartview=FitH, hyperfigures=false,
            plainpages=false, naturalnames=true,             	      colorlinks=true,pagebackref=true,
pdfpagelabels]{hyperref}
  \usepackage[numbers,sort&compress]{natbib}
 
\usepackage{color,xcolor}

\usepackage[T1]{fontenc}

\newtheorem{theorem}{Theorem}[section]

\theoremstyle{definition}

\newtheorem{thm}{Theorem}[section]

\theoremstyle{definition}
\newtheorem{dfn}[theorem]{Definition}

\theoremstyle{remark}
\newtheorem{rmk}[theorem]{Remark}

\numberwithin{equation}{section}

\def\lh{\hbox to 15pt{\vbox{\vskip 6pt\hrule width 6.5pt height 1pt}
\kern -4.0pt\vrule height 8pt width 1pt\hfil}} 

\def\qed{\hbox{${\vcenter{\vbox{\hrule height 0.4pt\hbox{\vrule width
0.4pt height 6pt \kern5pt\vrule width 0.4pt}\hrule height 0.4pt}}}$}}

\def\C{\mathcal C}

\def\C2{F(S^1,2)}
\def\C3{\overbar F(S^1,3)}
\def\C4{\overbar F(S^1,4)}

\def\overbar{\overline}

\def\4mu4{\nu^{a}_{\alpha\beta}}

\def\Aoo{A_\infty}

\def\be{\begin{eqnarray}}
\def\ee{\end{eqnarray}}

\def\ddo{\end{document}}

\newcommand{\etc}{\emph{etc.}, }
\newcommand{\ie}{\emph{i.e.}, }
\newcommand{\eg}{\emph{e.g.}, }
\newcommand{\cf}{\emph{cf.}, }

\def\shreps{representations up to homotopy / homotopy coherent representations}

\begin{document}
\title{Representations up to coherent homotopy}
 \author{Tim Porter}
 \address{Ynys M\^{o}n /  Anglesey, Cymru / Wales, ex-University of Bangor}
  \email{t.porter.maths@gmail.com}
\author{Jim Stasheff}
\address{University of Pennsylvania}
\email{jds@math.upenn.edu}

\begin{abstract}
 \emph{Homotopy coherence} has a considerable history, albeit also by other names. 
We provide a brief semi-historical survey providing some links that may not be common knowledge.\footnote{\today}
\end{abstract}
\maketitle

\tableofcontents

\section*{\emph{Avant propos.}}

 What is a \emph{symmetry}? What is a representation?
 
 `\emph{The time has come\footnote{ the walrus said -  Lewis Carroll} to talk of many things}', which are really different aspects of a single set of ideas from Homotopy theory, Galois theory, Non-abelian Cohomology, Grothendieck's \emph{Pursuit of Stacks}, and more. We will look at how generalisations of the ideas of symmetry and of what  representations are, leads to a thread that weaves its way through these areas of mathematics.
 
In mathematics, symmetries are embodied in the action of a group, $G$, on a mathematical object,  $X$, so there is some $\alpha: G\times X\to X$, satisfying the usual rules; better yet, a homomorphism from $G$ to $Aut  X$, the group of invertible maps from $X$ to itself, and even better \ldots , but that will have to wait. The objects may be \eg just sets, geometric figures or vector spaces, etc., and we \emph{represent} the elements of the group by permutations, or matrices, or \ldots . On the algebraic side of  math-physics, the objects are often `spaces of fields' and one encounters \emph{infinitesimal} symmetries, referring to an action of a Lie algebra on the fields. In homotopy theory or homological algebra, instead of a homomorphism, one may have a map, $\alpha: G\to Aut X$, which is a homomorphism only up to some form of homotopy. For example,
though singular and de Rham cohomology agree on smooth manifolds, the comparison at the cochain level is multiplicative only up to homotopy; in fact, there are further  higher homotopies as developed by Gugenheim and Munkholm.
The compatibility of these higher homotopies is now expressed as \emph{homotopy coherence}.


\section{Introduction}

To take just one such occurrence of symmetries and representations arising in a topological setting, one that we will explore a bit more fully in the next section, in the classical theory of covering spaces of a (pointed) space, $X$, one has the result that covering spaces of $X$ and actions of the fundamental group, $\pi_1(X)$, on discrete spaces correspond. If one thinks of covering spaces as fibrations with discrete fibres, then one may try to extend such an action to an action on the fibres of a more general fibration over $X$, but that does not work well because of all the choices that have to be made. If, however, one replaces $\pi_1(X)$ by the based loop space, $\Omega X$, then a form of action does remain. This idea can be extended, ultimately freeing the theory of the need for a base-point  and replacing $\Omega X$ by the singular complex of $X$.  We will touch on this later on.  Of course, $\Omega X$ is not a group, but it is in many respects a homotopy coherent form of a group, and the action that remains on the fibre over the base point is a homotopy coherent action. The algebraic operations can be defined by making choices, but those choices are linked by homotopies, and homotopies between homotopies, \etc  so that the fibration corresponds to a \emph{homotopy coherent representation} of $\Omega X$,  but what does that mean, and  what is the story  behind the term `homotopy coherence'?

To our knowledge, the first occurrence of homotopy coherence (without the name) was in the work of Sugawara, \cite{sugawara:g,sugawara:h, sugawara:hc}. Even in one of the early source papers, Vogt, \cite{vogt:1973}, the wording was ``up to coherent homotopies''. The earliest use of ``homotopy coherence'', as such and
that we know of, is due to Cordier, \cite{cordier82}. The study initiated by Sugawara in the context of topological $H$-spaces grew over time along parallel threads and recently has expressed new vigour in relation to \emph{$\infty$-local systems.}

In this survey article, we will sketch out various parts of the theory of homotopy coherence and its application to  ideas of representation of a homotopy type and the related idea of an $\infty$-local system. This will continue, in part, some of the threads that were mentioned in our earlier paper, \cite{TP-JDS-hcr:2022}, but with some areas being treated with a different emphasis and others not followed up at all.
 
 Some of these areas will be in a topological setting, but others in the closely related area of homotopical algebra and in the applications of these ideas to differential geometry and to some limited  extent to theoretical physics. We will also explore the covering space example, as it leads on to the higher dimensional homotopy coherent representations  being linked to certain types of `spaces over a base' and thus to stacks, gerbes, and Grothendieck's \emph{Pursuing Stacks}, \cite{groth:PS}, although we will not comment  on this in much detail, due to lack of space and time.
 
This survey will be  semi-historical  and idiosyncratic  with the topics covered determined by the knowledge and taste of the authors, but we hope it will provide some links that may not be common knowledge between the various aspects of the theory of homotopy coherence and, in particular, to its application in the form of homotopy coherent representation theory.

If one examines this subject area in any depth, one realizes that there are various themes that intertwine in the historical development of the theory. There is a topological theme, but, in the development of that, one naturally finds a  subtheme that passes through the land of  chain complexes, as applied via  (co)homological methods to problems in topology and algebra. This, then, develops  a life of its own in more purely algebraic areas via the topics of differential graded (dg) categories and later $A_\infty$-categories and algebras. Leading out from the topological theory,  there is also a simpliciaily based subtheme that  eventually leads to the study of simplicially enriched categories, quasi-categories and various manifestations of $\infty$-category theory.  To some extent, these various  themes have developed semi-independently.
It seems inevitable, therefore, that 
this may result in a somewhat disjointed flow of concepts as the developments are to some extent  parallel, but that \emph{does} reflect the history of the subject.  Other parts of the story involve input from algebraic geometry and non-abelian cohomology with the Grothendieck construction and, more recently, some parts of Grothendieck's Pursuit of Stacks and, consequentionally, input from (higher) category theory. 

All these threads  have, however, recently started to converge in some of their methods. We will attend to aspects of this convergence by working towards a description of homotopy coherent / $A_\infty$-local systems on a space and possible homotopy coherent generalization of the classical Riemann-Hilbert theory.

As this is a survey, we do not claim to be saying anything really new. We do, however, think that as the theory, or more exactly theories, have pushed further forward, it is important to be reminded of links, questions  and methods that were apparent in earlier parts of the history, but which have to some extent been sidelined, yet can contribute a useful perspective on more modern results.

\section{Background/Motivation}
There are at least three different directions from which to approach the basic ideas that we will be exploring.
\subsection{The (classical) topological setting}\label{Background top}
In topology, the fundamental theorem of covering spaces asserts that for a `nice'\footnote{globally and locally path-connected and semi-locally simply-connected} (pointed)  topological space $X$,  the functor\footnote{Way back when, functors had not yet been named!} given by sending a covering space over $X$ to the corresponding representation of its fundamental group $\pi_1(X)$ as permutations of a given discrete fiber over a point is an equivalence.
That is, the representation determines 
 a covering space for which the monodromy/holonomy
 is naturally isomorphic to the original representation.

Let us recall that the theory of covering spaces provides the link between two classical definitions of the fundamental group of a space, $B$. In one, which may be called the \emph{Poincar\'{e} fundamental group} of $B$, we choose a base point, $b\in B$, and $\pi_1(B,b)$ will be the group of homotopy classes of paths from $b$ to itself.  Provided that $B$ is path connected, this does not depend, up to isomorphism, on which  point is chosen as base point.

 There are some obvious generalizations of this. We could replace a single base point by a set, $B_0$, of base-points, so getting a groupoid version, which would work even when the space is not connected.  We could replace the use of homotopy classes, working instead with the paths themselves, and this would give some version of the loop space, $\Omega(B,b)$, which has already been mentioned, and finally we could replace paths or loops by higher dimensional things, so getting something, perhaps the $n^{th}$ homotopy groups, $\pi_n(B,b)$, or some relative analogue. These however, for now, do not have any immediate analogue of the covering space theory associated to them. 
 
 The process of going from the base-pointed case to the multiple object case is one of the features of \emph{categorification}, but for \emph{higher categorification}, we need to replace sets by categories, categories by 2-categories, and so on.  Remember  that, in this process, strict higher categories are not sufficient to model the general case, and that leads to considerations of homotopy coherence, the modelling of homotopy $n$-types and aspects of $\infty$-category  theory. Of course,  the resulting theory should, hopefully, interpret and link with higher order geometric structures, both in the mathematics and in the related mathematical physics.

 The second version of the fundamental group is variously called the \emph{Chevalley fundamental group} or the algebraic fundamental group, and relates to what has been called the \emph{Galois theory of covering spaces.} For the moment, we will only consider it for connected spaces, $B$, which admit a universal covering map, $p:E\to B$, in which $E$ is connected.  This fundamental group is defined to be $\mathsf{Aut}(p)\equiv \mathsf{Aut}_B(E,p ),$ the group of automorphisms of the universal covering, $(E,p)$, of $B$.
 
 For many of the usual connected spaces that are met in homotopy theory or differential geometry, the two definitions give isomorphic groups, and, of course, this is a consequence of the classical theory of coverings.
 
 This `algebraic' version of the fundamental group is, thus, as the automorphism or symmetry group of an object in the category, $Cov(B)$, of coverings of $B$.  We can adjust that viewpoint usefully by introducing a base point, $b$, in $B$ and also  the fibre \emph{functor},
 $$Fib_b:Cov(B)\to Set,$$
 which to a covering, $q:C\to B$, assigns the fibre over $b$, \ie $Fib_b(C,q)=q^{-1}(b)$. An alternative, and equivalent, definition of the `algebraic' fundamental group of $(B,b)$ is that $\pi_1^{alg}(B,b)$ is the group of automorphisms of $Fib_b$. 
 
 Here we have a formulation that is more versatile for generalization.  We could pick a set of base points and look at the groupoid of invertible natural transformations between the corresponding fibre functors; that is easy to do.  The other potential generalisations that we mentioned for the Poincar\'e version need a bit more thought, but do work out provided we move away from $Cov(B)$ in which the fibres are sets /  discrete spaces towards more general fibred objects. This process shows some of the inadequacies of any theory considering just the higher homotopy groups on their own and not the detailed homotopy structures that those groups encode. The higher homotopy groups do not tell us enough about the homotopy types involved,  needing the encoding of extra homotopical information. A map between CW-complexes which induces isomorphisms on all homotopy groups  is, by the famous theorem of J. H. C. Whitehead, a homotopy equivalence, but the existence of isomorphisms does not guarantee the existence of such a continuous mapping; see discussions of this in most textbooks on homotopy theory.   The covering space theory involves \emph{representations} of $\pi_1(B)$ on sets, so of the homotopy 1-type of $B$ on families of sets.  To go further, and to model, say, the homotopy 2-type by representations up to homotopy requires a bit more machinery. 
 \begin{rmk}
 Recall that a space or simplicial set is considered to be a \emph{homotopy $n$-type} if its homotopy groups vanish  above dimension $n$, so a group represents a homotopy 1-type, a crossed module represents a homotopy 2-type, and so on. \end{rmk}
 
\begin{rmk} We should also mention the way in which Grothendieck, in SGA1, \cite{SGA1}, adapted the Chevalley approach to handle the case of schemes, the type of space, endowed with a structural sheaf of rings, that arises naturally in Algebraic Geometry. There is an analogue of the theory of covering spaces in that context that replaces considering all coverings by just the \emph{finite} coverings of a space. There is not usually a universal covering map available, but Grothendieck adapted a categorical version of the Chevalley theory, (replacing representability by pro-representability), to arrive at a notion of a fibre functor and an automorphism group that has a natural profinite topology on it.  In certain cases, this profinite fundamental group is the Galois group of a field extension, and, in general, there is a continuous representation of the profinite fundamental group on
finite discrete spaces.  This theory of Grothendieck and its relationship with both Galois theory and with the theory of covering spaces is well explained in the books by Douady and Douady, \cite{DouadyR:algetg:1974}, Szamuely, \cite{Szamuely:Galois:2009}, and Borceaux and Janelidze, \cite{FB&GJ}. \end{rmk}

 To return to the main thread, but now going in the direction of higher homotopies and  both $\infty$- and $A_\infty$-category theory, see section \ref{A-infty}, we need to look at fibrations over $B$ with fibres that are more complicated than discrete one. At its most general, this line of enquiry leads to Grothendieck's \emph{Pursuit of Stacks} and to links with non-abelian cohomology, but long before we can hint at that, we need to go back and look again at fibrations in the classical theory.

We note that there are results similar to those of covering space theory, which link smooth fibre bundles with flat connection and representations on vector spaces,
 but for  (Hurewicz) fiber spaces\footnote{We will use the terms \emph{fibration} and \emph{fibre space} to indicate the presence of some   Homotopy Lifting Condition, which we will discuss more fully later on.}, there is a much more subtle correspondence.

 Consider such a fibre space, $F\to E\to B$, where $F$ is the fibre over a base point $b\in B$.
 One version of the classification of such bundles/fibrations is by the action of the based loop space, $\Omega B$, on the fiber $F$. This began with Hilton, \cite{hilton}, showing that there was an `action' $\Omega B \times F \to F$, which at the level of homotopy classes gave a representation of  $\Omega B$ on $F.$
 In fact, that action gives rise to two maps,  $\Omega B \times \Omega B \times F \to F$, which are homotopic.
  However, by passing to the homotopy classes of based loops in $B,$ that information is lost; there was more to be discovered. In those days,  
  $\Omega B$ was the space of based loops which were parameterized in terms of the unit interval, hence composition was associative only up to homotopy.
   It was Masahiro  Sugawara, \cite{sugawara:g,sugawara:hc}, who first studied `higher homotopies',  leading to today's world of $\infty$-structures of various sorts \cite{loday-vallette,lurie:HighTop2006,sss}  (see section \ref{A-infty}), and even a \emph{Journal of Higher Structures.}

Although Sugawara did not treat actions up to higher homotopies as such, he did introduce \emph{strongly homotopy multiplicative maps} of associative H-spaces. With John Moore's introduction of an associative space of based loops, the adjoint map, $\Omega B  \to F^F$, could be seen as such a map. Trying to consider the action as a \emph{representation} of $\Omega B$ on $F$, we have perhaps the first example of a \emph{homotopy coherent representation}. The algebraic structure on the map from $\Omega B$ to $F^F$ is `up to homotopy' and those homotopies fit together, again `up to homotopy'. Equations in that structure are replaced by explicit homotopies, themselves linked by higher homotopies. What is more, the `action' is similarly `up to explicit homotopies', and higher homotopies, see below and section \ref{homotopy coherent nerve}.

 \begin{rmk} An apparently very different sort of treatment of \emph{higher homotopies} occurred, slightly later,  in a strong form of Borsuk's shape theory.  As we have already explored this to some extent in \cite{TP-JDS-hcr:2022}, we will not repeat our discussion in detail, but make reference to Marde{\v{s}}ić's book, \cite{mardesic:strong:2000}, and to `the literature'.\end{rmk}

The link between that geometric topological approach and the notions emerging from Sugawara's work came from the work on \emph{homotopy everything algebraic structures} by Boardman  and Vogt,  \cite{boardmanvogt},  then in Vogt, \cite{vogt:1973}, which gave the connection between a detailed `geometric' approach to homotopy coherent diagrams and a homotopical approach related to model category theory.  We will see further links to Vogt's results later on.

 Associativity being central to Sugawara's work on homotopy multiplicativity, \cite{sugawara:hc}, it was not long before  maps of spaces with associative structures were  generalized to maps of topological categories and \emph{homotopy coherent functors}, which 
arise in other contexts involving topological  and simplicially enriched categories, \cite{C&P97}.

In topological settings, one often uses that the categories involved have enriched structures,  such as being topologically or simplicially enriched. We note that topologically enriched categories are sometimes called topological categories, although that term has other meanings as well, whilst simplicially enriched categories are sometimes called  simplicial categories, and, again, that is ambiguous. In these enriched cases, higher homotopies have good relatively simple interpretations.  The following is adapted from Vogt's description.

For  a category,  $\mathcal C$, in general, we will denote by  $\mathcal C_n$ the subset of  $\mathcal C^n$ consisting of composable strings of morphisms, \eg $C_1 = Mor\   \mathcal C$. We use juxtaposition to denote the composition.
\begin{dfn} For topologically enriched  categories, $\mathcal C$ and $\mathcal D$, a \emph{functor up to strong 
homotopy} $F$, also known as a \emph{ homotopy coherent functor}, consists of maps $F_0: Ob\ \mathcal C\to Ob\ \mathcal D$, $F_1:Mor\ \mathcal C \to Mor\ \mathcal D $ and, for each higher value of $n$,
maps $F_n:I^{n-1}\times {\mathcal C}_n \to {\mathcal D}_n$  such that 
  \begin{alignat*}{2}
  F_1(&x\to y) : F_0(x)\to F_0(y), \notag \\
  \intertext{and}
F_p(t_1, &\ldots ,t_{p-1}, c_1,\ldots
 ,c_p)\\& = \begin{cases}F_{p-1}(\ldots,\hat{t}_i,\ldots,c_ic_{i+1},\ldots
 ) &\mbox{if } t_i=0  \notag\\
F_i(t_1, \ldots ,t_{i-1}, c_1,\ldots,c_i)F_{p-i}(t_{i+1},\ldots ,t_{p-1},c_{i+1},\cdots ,c_p) &\mbox{if } t_i=1.  \notag
 \end{cases} 
 \end{alignat*}
  \end{dfn}
 `Homotopy coherence' thus relates to how the maps on the faces of these cubes fit together.
We will encounter  this idea of compatible maps on cubes later.
\subsection{... and for dg categories?}\label{Background dg}

There are complete analogues of these ideas for categories of (co)chain complexes and, more generally, differential graded (dg) categories.
We note that for  dg-categories (that is, categories enriched over some category of chain complexes)  a chain analogue of cubes applies, so the $I^{n-1}$ in the above is replaced by a chain complex  representing it.

This allows one to transfer, or translate, certain problems from the topological setting to that of (co)chain complexes, but these notions have been found to be relevant for many other contexts without that immediate topological link, in particular within representation theory and homological algebra.


For an associative algebra $A$, an $A$-module, $M$, can be considered as a representation of $A$, that is a morphism, 
$A\otimes M\to M$, or $A\to End(M)$.  In a differential graded context, one can consider representations up to homotopy of $A$ on $M$.
On the chain (dg) level, the corresponding notion is related to that of a \emph{twisting cochain}, the twisted tensor product differential as
 introduced in \cite{brown59}, for
modelling the chains on the total space of a principal fibration in terms of chains on the base and chains on the fibre, (see  section \ref{twisting cochains}).

One can up the ante further by introducing $A_\infty$-spaces, $\Aoo$-algebras and modules and further going towards $\Aoo$-categories (see section \ref{A-infty}), 
but take care: do you want objects and morphisms to be $\Aoo$ in an appropriate sense or only the morphisms?  The most general case is not always the most useful one for applications, as the level of complication in some of the formulae can lead to `diminishing returns'!

 It is worth pointing out that $A_\infty$-maps between strictly associative dg-algebras were studied earlier than  \cite{jds:hahI,jds:hahII},   by Sugawara, as being  parameterized by cubes.  They are the  \emph{strong homotopy muliplicative} maps that were mentioned earlier.  $A_\infty$-maps of  
 $A_\infty$-spaces are parameterized by more complicated polyhedra.  We will return to some of these ideas a bit later on in section \ref{A-infty}.
\subsection{Fibred categories and pseudo-functors}
Further steps towards the higher categorification of these ideas came from Algebraic Geometry and non-abelian cohomology. This had started independently of the above, back in the late 1950s.

In the Galois-Poincaré theory, or perhaps it should be the  Galois-Chevalley-Grothendieck theory, of covering spaces and actions,  the  covering spaces are fibrations, but with discrete fibres, \ie the fibres are sets, and, from the point of view of homotopy theory, sets only model the information on the connected components of a space, that encoded by $\pi_0$.  

A first generalisation of this is to look at fibrations whose fibres are (homotopy) 1-types, thus specifying both $\pi_0$ and $\pi_1$, or better still, specifying the fundamental groupoid.  A covering space is given by a locally constant sheaf of sets, so a generalisation might be a sort of categorification, or rather groupoidification, of a sheaf of groupoids.  Such a notion occurs naturally in Algebraic Geometry, where Grothendieck, again in SGA1, \cite{SGA1}, but also earlier in \cite{Grot:TechDescI}, identified the notion of a fibred category, and, more particularly, a category fibred in groupoids.  These are given by \emph{pseudo-functors}. Recall that these behave almost like functors, but preserve composition and identities only up to equivalence.  These pseudo-functors will take values in the 2-category, $\mathsf{Grpd}$, of groupoids, functors between them and natural equivalences between those.
The 2-category has `hom-sets' which are groupoids. Later on we will see that these pseudo-functors are precursors of the general notion of a representation up to homotopy\footnote{as the homotopy structure of the category of groupoids is based on natural transformations.}. The intended intuition is that a fibred category should be some sort of fibration of categories, so should involve a family of categories indexed by the objects of a `base' category, and related by functors `over' the morphisms of that base.\label{non-abelian}

To link up, more fully, pseudo-functors with the topological and simplicial parts of the story, we will need to recall the nerve construction.  We denote the ordered set $\{0<1<\ldots, <n\}$ by $[n]$, and if $\mathsf{C}$ is a small category or groupoid, then $Ner(\mathsf{C})$ is the simplicial set having as its $n$-simplices the functors from $[n]$ to $\mathsf{C}$,  with face and degeneracy maps induced in the obvious way; see almost any source on simplicial homotopy theory, but, in particular, Joyal's notes, \cite{joyal:2008}, as that source also handles the special properties enjoyed by the nerves of categories. We will use that $Ner$ is a functor from the category of small categories to that of simplicial sets.

 By applying the nerve functor to each hom-groupoid, one obtains that $\mathsf{Grpd}$ can also be considered as a simplicially enriched category having Kan hom-objects, the most important condition for a simple description of homotopy coherence, \cf for instance,  \cite{cordierporter:Vogt:1986}. The main point to note is that these categories fibred in groupoids are what might be called `2-local systems' or `2-representations' that is, representations up to natural equivalence. 

\begin{rmk}
A good introductory source for the theory of categories fibred in groupoids and the related notions of 1-stacks, is Breen's notes, \cite{breen:2006}, of his course at the IMA session on higher categories in 2004.  We may refer to these notes several times as they will obviate the need to give too lengthy an exposition here.
\end{rmk}
The two notions of fibred categories / fibration of categories and pseudo-functors are more-or-less equivalent. A Grothendieck fibration (of categories) is a functor, $p:\mathcal{E}\to \mathcal{B}$, having certain lifting properties that we will mention later.  This leads to a pseudo-functor from $\mathcal{B}$ (or its opposite) to $\mathsf{Cat}$, the 2-category of small categories.  The fibration is \emph{fibred in groupoids} if, for each object, $b$, of $\mathcal{B}$, the fibre category, $p^{-1}(b)$, is a groupoid, or, equivalently, if the pseudo-functor factors through the inclusion of $\mathsf{Grpd}$ into $\mathsf{Cat}$.

The pseudo-functor approach is a generalisation of the notion of a presheaf of local sections of a space over a base, which is one form of a (fairly classical)  local system. This  was, in fact, the original form of category fibred in groupoids given by Grothendieck in \cite{Grot:TechDescI}.  It is also the form used in Breen's notes, so it is convenient to follow that approach here. The two forms are linked by the Grothendieck construction  that will be mentioned more fully  on page  \pageref{Groth. constr.}.

We will concentrate on the case in which the base category, $\mathcal{B}$, is the category of open sets of a space,\footnote{As Breen notes, `space' here could be interpreted in many useful ways, \eg as a scheme, or as a site, that is, a small category together with a Grothendieck topology.} $X$.

\begin{dfn} (Pseudo-functor version)
A category fibred in groupoids over a space, $X$, is given by a family, $\{\mathcal{C}_U\}$, of groupoids indexed by the open sets, $U$, of $X$, together with 
\begin{itemize}\item an inverse image functor, $$f^*:\mathcal{C}_U\to \mathcal{C}_V,$$
for each inclusion, $f:V\to U$, of open sets (which is assumed to be the identity on $\mathcal{C}_U$ when $f$ is the identity on $U$);
\item for each pair of composable inclusions,
$$U_2\xrightarrow{g}U_1\xrightarrow{f} U_0,$$
a natural isomorphism,
$\varphi_{f,g}: (fg)^*\Rightarrow g^*f^*;$
\item[]\hspace{-11mm}and
\item (cocycle condition) for each sequence
$$U_3\xrightarrow{h}U_2\xrightarrow{g}U_1\xrightarrow{f} U_0,$$ the composite natural isomorphisms,
$$(fgh)^*\Rightarrow h^*(fg)^*\Rightarrow h^*(g^*f^*),$$
and
$$(fgh)^*\Rightarrow (gh)^*f^*\Rightarrow (h^*g^*)f^*),$$coincide.
\end{itemize}\end{dfn}
This is thus a pseudo-version of a presheaf, $\mathcal{C}: Open(X)^{op}\to  \mathsf{Grpd}$, of groupoids on $X$.  If the fibred category concerned also satisfies that arrows glue locally and objects glue, then it is a \emph{stack of groupoids}, and again see Breen, \cite{breen:2006}, for more precision.

Remembering that pseudo-functors preserve composition up to natural equivalence, and that, in $\mathsf{Grpd}$, this translates into the natural notion of homotopy, we have a version of `representation up to (coherent) homotopy'\footnote{As $\mathsf{Grpd}$  has only a low dimensional version of homotopy, the coherence is manifested in just the cocycle condition.}.
We can  easily give some examples. One of the usual ones, again given by Breen, \cite{breen:2006}, is that of the category of $G$-torsors / principal $G$-bundles specified locally on $X$, so $\mathcal{C}_U$ is the category of principal $G_U$-bundles on $U$, where $G$ is typically a sheaf of groups. 

If we go to the generalisation of pseudo-functor on an arbitrary category, then examples are provided by the Schreier theory of (non-abelian) group extensions, \emph{cf.} Blanco, Bullejos and Faro,  \cite{BBF:MZ:2005}, and by Haefliger's theory of complexes of groups, see Fiore, L\"{u}ck and Sauer, \cite{Fiore-Luck-Sauer:euler:2011}, which thus links with the theory of orbifolds.

Any functor will give a pseudo-functor, so for a groupoid, $\mathcal{G}$, and an action of it on sets, we have another example, thus any local system gives an example. Because $Set$ does not have a 2-category structure, any pseudo-functor to $Set$ is just a functor.

The other forms of the definition of fibred category are as a \emph{Grothendieck fibration} or \emph{op-fibration}. This is a functor, $p:\mathcal{E}\to \mathcal{B}$, with a \emph{lifting property} for arrows in $\mathcal{B}$, so specifying an object over the codomain of an arrow in $\mathcal{B}$, there is an arrow\footnote{Cartesian or op-cartesian arrows} in $\mathcal{E}$ over it, with nice universal properties. One of the classic examples is that in which $\mathcal{E}$ is the category of modules (over all rings), so with objects pairs, $(R,M)$, where $R$ is a ring and $M$ is a (left) $R$-module, $\mathcal{B}$ is the category of rings, and $p$ sends $(R,M)$ to $R$. If we have a ring homomorphism, $f:R\to S$, and an $S$-module, $N$, then one has $f^*(N)$ is $N$ with the induced $R$-module structure and $(R,f^*(N))\to (S,N)$ is a cartesian lift of $f$. This is a very rich example, but we leave the exploration to the reader.  In this example, we also have a pseudo-functor. In general, given a pseudo-functor, $F:\mathcal{B}\to \mathsf{Cat}$, one can build a Grothendieck fibration using a form of the \emph{Grothendieck construction}.

 Looking at the analog of local sections in this categorical setting, one finds that a fibred category, as a functor, $p:\mathcal{E}\to \mathcal{B}$, corresponds to a pseudo-functor
 from the base category to the category of small categories, picking out the fibres and the functors linking them together.

 The construction of the covering space from the action is a simple example of a type of semi-direct product.  Both Ehresmann and Grothendieck looked at the relevance of this for fiber spaces and for the fibred categories that Grothendieck had introduced in his work.  The notion of a fibred category mimics both topological fibrations and the relationship between the category of all modules (over all rings) and the category of rings itself. \label{Groth. constr.}
 Grothendieck's construction started with a pseudo-functor and constructed a fibred category. In other words, it started with a homotopy coherent representation of the base category and built a fibred category (aka Grothendieck fibration)  from it.  It is of interest that, by Thomason's results, \cite{thomason:thesis:1977,Thomason:1979}, the Grothendieck construction is a form of \emph{homotopy colimit} within the context of the natural homotopy theoretic structure of $\mathsf{Cat}$. 
 
In his extensive discussion document, \cite{groth:PS}, Grothendieck started exploring the `higher homotopy' analogues of fibred categories, thus analogues of simplicial fibrations.  The fibres were to be models of homotopy types, extending both  covering spaces,  categories fibred in groupoids, and, more generally,  homotopy $n$-groupoids, for any value of $n$ including $\infty$. These were his $n$-stacks, which can be thought of a homotopy coherent versions of (pre)sheaves of  homotopy $n$-types.

As it arose in this context,  it is appropriate to mention here Grothendieck's so-called\footnote{In fact, Grothendieck did not form such as a  \emph{hypothesis}, but rather put forward the idea  that a good test for any possible  model for homotopy $\infty$-groupoids should be that the objects should model all homotopy types, and that the higher categorical structures should also be equivalent. The term was coined considerably later (2007) and unfortunately seems to have taken root!} \emph{Homotopy Hypothesis}, namely that whatever meaning is given to  the term \emph{homotopy $n$-groupoid}, as above, it must be general enough to model all homotopy $n$-types.  This leads to the point made by one of us (TP) in a letter, \cite{Porter:Grothendieck-letter:1983}, to Grothendieck (16 June 1983) that Kan complexes fit the bill to be $\infty$-groupoids in this sense.  This is discussed in \cite{porter_:Spaces-as-infty1grpds-2021}.  If Kan complexes play that role in the theory, then to replace general simplicially enriched categories  by ones that are `locally Kan' makes a lot of sense, as, if Kan complexes are thought of as $\infty$-groupoids, then locally Kan  enriched categories look like $\infty$-categories in which all higher morphisms are weakly invertible.  Applying  Cordier's  homotopy coherent nerve functor then makes the use of weak Kan complexes / quasi-categories very natural, linking up the modern Joyal-Lurie theory of $\infty$-categories with homotopy coherence  in the Vogt-Cordier theory.

This makes it clear that one of Grothendieck's motivations for  \emph{Pursuing Stacks} was to give a generalisation to all levels of the classical correspondence between covering spaces and representations of the fundamental group of a space.   Homotopy coherence  of the representation of the base homotopy type is key to building the fibred structures over that base. 

\section{From $A_\infty$-spaces to $A_\infty$-categories}\label{A-infty}

So far in our quest for the interpretations and  meanings  of `representations up to coherent homotopy', we have assumed mostly that the object being represented was somehow algebraic, so was a group, groupoid or category, but we did point out early on that some of the objects naturally involved in `actions' were themselves `algebraic up to coherent homotopy' and it is to this we now turn as it represents a different parallel thread of the story.

A well known example of homotopy coherence (though not originally named as such) is the  structure on the loop space, $\Omega X$, of a based space, $X$. This space and its structure, of course, has many interpretations, but here we give just one as a fairly simple example of the rich structure we study further. 

There are multiple composition mappings, 
$$\lambda_n:(\Omega X)^n\to \Omega X,$$
given by  generalizations of  the usual one related to associativity up to homotopy:

$(ab)c$ is only  homotopic to $a(bc)$, 
via  an \emph{explicit} (chosen) homotopy, $\mu_3: I \times (\Omega X)^3 \to \Omega X$; we will look at higher associativity as well, but there's more going on here  as abstraction reveals how unbiquitous this sort of structure is.

Again the story starts in the 1950s, as was  mentioned in section \ref{Background top}.

The history of $A_\infty$-structures begins, implicitly, in 1957 with the work of Sugawara, \cite{sugawara:h,sugawara:g}.
He obtains criteria for  a space, $Y$, to be a homotopy associative H-space or a loop space.

When one of us, (JDS), was a graduate student at Princeton, his supervisor John Moore suggested he look at when a primitive cohomology class, $u\in H^n(Y,\pi)$, of a loop space, $Y=\Omega X$, was the suspension of a class in $H^{n+1}(X,\pi)$.  In other words, when was an H-map $Y\to K(\pi,n)$ induced as the loops on a map to $K(\pi,n+1)$.  JDS attacked the problem by considering Sugawara's 
infinite sequence of conditions, ``higher homotopies'', generalizations of homotopy associativity.  Sugawara also had included conditions involving homotopy inverses, which could be  avoided via mild restrictions on $Y$.

To systematize the \emph{specific} higher homotopies, $A_n$- and $A_\infty$-spaces  were introduced in \cite{jds:hahI}:

\begin{dfn}
An $A_n$-\emph{space}, $X$, consists of a space, $X$, together with a coherent set of maps, 
$$m_k: K_k \times X^k\to X\  \hbox { \it for  }\  2\leq k\leq n,$$
where $K_k$ is the (by now) well known $(k-2)$-dimensional associahedron.
\end{dfn}
The \emph{``so-called''} Stasheff polytope, $K_n$, was in fact constructed by Tamari in 1951, \cite{tamari:monoides:BSMF:1954,tamari:1962,huguet-tamari:1978}, a full decade before Stasheff's.

Many other realizations are now popular and have been collected by Stefan Forcey, \cite{forcey:zoo}, along with other relevant ``hedra''.
Tamari's point of view was much different from JDS's, but just as inspiring for later work, primarily in combinatorics. The book, \cite{Muller-Pallo-Stasheff:associahedra:2012}, has a wealth of offspring.

These higher homotopies were later referred to as a coherent family of  maps. In fact, they can now be said to form a homotopy coherent diagram of spaces, $(\Omega X)^n$, indexed by $n\geq 0$, and maps $\lambda_{r_1}\times \ldots \lambda_{r_n}:(\Omega X)^m\to (\Omega X)^n$ with $m=r_1+\ldots r_n$.
This gives a homotopy coherent diagram by using variants of the usual associativity homotopies for composition of paths.

Fukaya, in his study of Floer homology, \cite{Fukaya1}, was lead to formalize $A_\infty$-categories as a slight modification of $A_\infty$-algebras. Just as  a group can be considered as a one object category leading to the idea of a groupoid as a multi-object category, so an $A_\infty$-algebra should be considered as a one object $A_\infty$-category leading to its many object analogue in the idea of  an  $A_\infty$-category.  A short description of what  an $A_\infty$-category is thus as follows.

\begin{dfn}An \emph{$A_\infty$-category}, $\mathsf{A}$, with object set, $\mathbb{O}:= Ob(\mathsf{A})$, consists of a family, $\{\mathsf{A}(S,T): S,T\in \mathbb{O}\}$, of graded $k$-modules (usually vector spaces), together with unit elements, $1\in \mathsf{A}(S,S)$, for $S\in \mathbb{O}$, and maps,
$$m_n:\mathsf{A}(S_{n-1},S_n)\otimes \ldots \otimes \mathsf{A}(S_0,S_1)\to \mathsf{A}(S_0,S_n),$$
of degree $n-2$, for $n\geq 1$ such that the analogues of the conditions for the definition of (strictly unital) $A_\infty$-algebra are satisfied.\end{dfn}\label{A-infty cat}

Fukaya points out that the original motivation for  $A_\infty$-algebras was (the chains on) the space of based loops 
$\Omega X$, so, of course, the path space $PX = X^I$ can be regarded as leading to an $A_\infty$-category.
\section{Segal categories}
Segal categories are another instance of algebraic structures that are  up to coherent homotopy. They also provide one of the models for $\infty$-categories. They are bisimplicial sets with some extra structure.   
\subsection{The coherence problem and Segal's trick to get around it}
The reason that Segal categories arise as they do is best sought in the paper, \cite{segal:1974}, by Segal, although it is not there, but rather in \cite{DKS:1989} that they were introduced, yet not named as such. (In fact their first naming seems to be in Simpson's \cite{Simpson-Effective:ArXiv:1997}.) In \cite{segal:1974}, one of the main aims was to get `up-to-homotopy' models for algebraic structures so as to be able to iterate classifying space constructions, to form spectra 
for studying corresponding cohomology theories and to help `delooping' spaces where appropriate. Various approaches had been tried, notably that of Boardman and Vogt, \cite{boardmanvogt}.  In each case the idea was to mirror the homotopy coherent algebraic structures that occurred with loop spaces,  \emph{etc}.

As an example of the problem, Segal mentions the following.  Suppose $\mathcal{C}$ is a category and that coproducts exist in $\mathcal{C}$.  How is this reflected in the nerve of $\mathcal{C}$?  It very nearly acquires a composition law, since from $X_1$ and $X_2$, one gets $X_{12} = X_1\sqcup X_2$, and two 1-simplices,
$$X_1\rightarrow X_{12} \leftarrow X_2,$$
but $X_{12}$ is only determined up to isomorphism.
Let $\mathcal{C}_2$ be the category of such diagrams, \ie  in which the middle is the coproduct of the ends.  There is a functor,
$$\delta_2: \mathcal{C}_2 \to \mathcal{C}\times \mathcal{C}$$and this is an equivalence of categories, but there is also a `composition law', $$m : \mathcal{C}_2 \to \mathcal{C}$$ given by picking out the coproduct. This looks fine, but, in fact, this tentative multiplication again hits the problem of associativity.   Segal's neat solution was to side-step the issue.  He formed a category, $\mathcal{C}_3$,  consisting of all diagrams of form,
$$\xymatrixrowsep{2ex}\xymatrixcolsep{2ex}
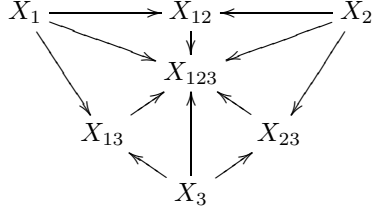
\xymatrix{X_1\ar[rr]\ar[drr]\ar[ddr]&&X_{12}\ar[d]&&X_2\ar[ll]\ar[dll]\ar[ddl]\\
&& X_{123}&&\\
&X_{13}\ar[ur]&&X_{23}\ar[ul]&\\
&&X_3\ar[ul]\ar[uu]\ar[ur]&&}$$
the notation being to indicate that each split line corresponds to the middle term being the coproduct of the two ends. All the usual natural isomorphisms between multiple coproducts are encoded in the one category.  There is an equivalence of categories 
$$\delta_3: \mathcal{C}_3 \to \mathcal{C}\times \mathcal{C}\times \mathcal{C}$$
sending the above diagram to $(X_1,X_2,X_3)$ and a `ternary operation' $\mathcal{C}_3 \to \mathcal{C}$ sending the diagram to $X_{123}$ compatibly, up to specifiable homotopies, with the structure outlined earlier.  The advantage is that all of this can be encoded by the nerve, and thus by the classifying space structure as a $\mathbf{\Gamma}$-space. 
\subsection{$\mathbf{\Gamma}$-spaces, $\mathbf{\Gamma}$-categories.}

What are $\Gamma$-spaces?
\begin{dfn}
(i) The category $\mathbf{\Gamma}$ is the category whose objects are all finite sets and whose morphisms from $S$ to $T$ are the maps $\theta :S \to \mathcal{P}(T)$ such that when $a\neq b\in S$, then $\theta(a)\cap \theta(b) = \emptyset$.  The composite of $\theta: S\to \mathcal{P}(T)$ and $\phi :T \to \mathcal{P}(U)$ is $\psi : S \to \mathcal{P}(U)$, where $\psi(a) = \bigcup_{b\in\theta(A)}\phi(b)$. 

(ii) A \emph{$\mathbf{\Gamma}$-space} is a functor $A: \mathbf{\Gamma}^{op}\to Top$ such that\\
(a) $A(0)$ is contractible;\\
and \\
(b) for any $n$, the map $p_n : A(n) \to\underbrace{ A(1)\times \cdots \times A(1)}_n$, induced by the maps $i_k : 1\to n$ in $\mathbf{\Gamma}$ where $i_k(1) = \{k\}\subset n$, is a homotopy equivalence.\end{dfn}

Our task, here, is not to review the contents of Segal's 1974 paper, so we will not explore $\Gamma$-spaces more  but rather to  note  the definition of a $\mathbf{\Gamma}$-category, which  follows the same model:

A $\mathbf{\Gamma}$-category is a functor, $\mathcal{C}: \mathbf{\Gamma}^{op}\to Cat$, such that (a) $\mathcal{C}(0)$ is equivalent to a one arrow category, and (b) as before except weak homotopy equivalence is replaced by equivalence of categories.
\subsection{Segal categories themselves}

The theory of $\mathcal{S}$-categories is sometimes too strict for convenience.  The notion of Segal 1-category is  a weakened version of that of $\mathcal{S}$-category.   The following is adapted from Toen's \cite{toen:cahiers:2002}, see also Simpson's text, \cite{simpson:homotopy-book:2011}.

The structure is first given, in their terminology, as a Segal 1-precategory, whilst Segal categories are precategories  satisfying the \emph{Segal condition}.
\begin{dfn}
\begin{itemize}
\item  A \emph{Segal 1-precategory} is specified by a functor
$$A : \mathbf{\Delta}^{op} \to \mathcal{S}$$
(\ie a bisimplicial set) such that $A_0$ is a constant simplicial set called \emph{the set of objects of} $A$.
\item A \emph{morphism} between two Segal 1-precategories  is a natural transformation between the functors from $ \mathbf{\Delta}^{op}$ to $\mathcal{S}$.
\item A Segal 1-precategory, $A$, is a \emph{Segal category} if for each $[p]$, the Segal map
$$\delta[p]: A_p \to A_1\times_{A_0}A_1  \times_{A_0}\ldots \times_{A_0}A_1,$$
is a  homotopy equivalence of simplicial sets.
\end{itemize}\end{dfn}
Segal categories form, together with quasi-categories\footnote{for which see  section \ref{qcat}.}, two of the main forms of model of $\infty$-categories, and as such correspond to another, equivalent,  interpretation of homotopy coherence.

\section{Homotopy twisting cochains}\label{twisting cochains}

Returning briefly to the subject of  fibered categories,  if we apply the nerve functor to a fibred category,  we will get something like a Kan fibration, so we should expect that similar ideas should be relevant when looking at simplicial fibrations and fibre bundles and, of course, they are.

The analog of a fibre bundle for simplicial sets,  $Y\to E\to B$, was analysed by  Barratt, Gugenheim and Moore, \cite{BGM}, using a \emph{twisting function}, $t:B\to aut(Y)$, and leading to a \emph{twisted cartesian product}.  
 A bit more precisely, given a base simplicial set, $B$, and a simplicial group, $G$, acting on a fibre $Y$, there is a simplicial set, $E$, that looks almost like the cartesian product, $B\times Y$, except that the $0^{th}$ face  is twisted by  twisting function, $t:B\to G$, so $d_0(b,y) = (d_0b, t(b)(d_0y))$. In fact, $t$ is given by a simplicial map, $t:B\to \overline{W}(G)$, for which see standard texts on simplicial homotopy theory.  If $Y$ is a discrete simplicial set, so really a set, this is a simplicial covering space, and the theory is, once again, that of a representation.

For a given simplicial group, $G$, there is a universal twisting function, which corresponds to a universal principal $G$-bundle, $W(G)\to \overline{W}(G)$, for which see standard sources on simplicial homotopy theory, such as Curtis, \cite{Curtis71} or May, \cite{May:1967}.

We will  next see what this looks like from the perspective of their associated chain complexes.

 In 1959,  Ed Brown, \cite{brown59}, constructed an 
 algebraic chain model of a twisted cartesian product, 
  a \emph{twisted tensor product}  of chain complexes, $(C_*B\otimes C_*F,D)$,
where  the usual untwisted differential, $d_B \otimes 1 + 1 \otimes d_F$, of a tensor product is twisted by adding a \emph{twisting term}, $\tau: C_*B\to C_{*-1} Aut(F)$, giving
$$
D = d_B \otimes 1 + 1 \otimes d_F + \tau
$$
 where $\tau$ corresponds to  a representation of
  the cobar algebra, $\Omega C_*B$, on $C_*F.$ 
  Here $C_*B$ is regarded as a dg-coalgebra and $C_* Aut F$ as a dg-algebra.
 The defining relation that $\tau$ satisfies is
 $$d_F \tau + \tau d_B = \tau \cup \tau.$$
 In turn, Szczarba, \cite{szczarba61}, showed that for a twisted cartesian product with  base $B,$ and  group $G$, 
there is a twisting cochain in $Hom(C_*B, C_*G).$ Kadeishvili, \cite{kad:twisting}, studied twisting \emph{elements} in relation to $a\smile_1 a$ in a  homotopy Gerstenhaber algebra,
 \cite{GerVor95}. Quite recently, for $G$-bundles, Franz, \cite{franz:szcz}, proved that the map $\Omega C_*B \to C_* G$ 
is a quasi-isomorphism of dg bi-algebras. That $\Omega C_*B$ is a dg bi-algebra follows from work 
 of Baues, \cite{baues:doublebar}, and Gerstenhaber-Voronov,  \cite{GerVor95}, using the homotopy Gerstenhaber structure of $C_*B$.

 For an algebra, $A$, in general, there is a \emph{universal} twisting  morphism, $\tau_A:BA\to A$,
such that for any twisting morphism, $\tau:C \to A$,  there is a unique
morphism $f_\tau:C \to BA$ of coalgebras with $\tau=\tau_A\circ f_\tau$, see \cite{HMS74}.

The above can be generalized to  \shreps\footnote{If we write $D\tau$ for $d_F \tau + \tau d_B,$ we can see this as another manifestation of the Maurer-Cartan principle,
 which has subsumed the integrablity condition in deformation theory and elsewhere.}, leading to a notion of a \emph{homotopy twisting cochain} with values in an $A_\infty$-algebra.
  
 \section{Cordier's Homotopy Coherent Nerve and infinity local systems}\label{homotopy coherent nerve}

The description of homotopy coherent diagrams, in probably the most standard form, is for a simplicial enriched category, $\mathcal{C}$, in which the hom-sets between the objects are Kan complexes. This fits the cases of $\mathcal{C}=Top$ or $Kan$ and can also  be adapted to handle the category of chain complexes with a bit of extra work.  The basic construction is due to Cordier, \cite{cordier82}, using ideas developed by Boardman and Vogt, \cite{boardmanvogt}.
\subsection{Cordier's nerve (1982)}
The idea is that for each finite ordinal,  $[n]=\{0<1< ... <n\}$, one builds a simplicial enriched category, $S[n]$.  This, in a certain sense, resolves $[n]$ regarded as a category.  In fact, this is almost just the comonadic resolution for the neat comonad forming free categories on underlying graphs.  Applied to $[n]$, we have $S[n]$ is the simplicially enriched category having the simplicial set  of maps from $i$ to $j$ being an $(j-i-1)$ cube, $\Delta[1]^{j-i-1}$. This gives that  $\mathcal{S}\!-\!Cat(S[n], \mathcal{C})$ gives the collection of all homotopy coherent diagrams within $\mathcal{C}$ having the form of an $n$-simplex. These collections form a simplicial class. (They will often be too large to be, strictly speaking, a simplicial set). This simplicial class is    called the \emph{homotopy coherent nerve} of the $\mathcal{S}$-category $\mathcal{C}$ and will be denoted $Ner_{h.c.}(\mathcal{C})$. It then is clear that a homotopy coherent diagram having the  form of a small category, $\mathsf{A}$, is given by a simplicial map from $Ner(\mathsf{A})$ to $Ner_{h.c.}(\mathcal{C})$. The construction of $S[n]$ can be extended to all small categories. An equivalent way of viewing a homotopy coherent diagram is as an enriched functor from $S(\mathsf{A})$ to $\mathcal{C}$.

\begin{rmk} (On terminology.) Cordier, \cite{cordier82},  used the term `homotopy coherent nerve' for the above as he was  interested in its use in the area of strong shape theory and related areas. In his subsequent work with TP (in Cordier-Porter, \cite{C&P88,C&P:1990,C&P96,C&P97}), the quasi-categorical and $\infty$-categorical aspects became more evident.  Lurie, \cite{lurie:HighTop2006}, has called this form of nerve the \emph{simplicial nerve functor} as his applications are not explicitly concerned with homotopy coherence.\index{simplicial nerve of an $\mathcal{S}$-category}  Various other names have been given to it in the literature, sometimes `reinventing' the notion. We will stick with the original term as the others that are being used neglect the homotopy coherent aspect.\end{rmk}

The history of this construction is interesting. A variant of it, but with topologically enriched categories as the end result, is in the work of Boardman and Vogt, \cite{boardmanvogt}, and also in Vogt's paper, \cite{vogt:1973}.  Segal's student, Leitch used a similar construction to describe a homotopy commutative
 cube (actually a \emph{homotopy coherent cube}), \cf \cite{leitch}, and this was used by Segal, \cite{segal:1974}, under the name of the `explosion' of $\mathcal{C}$. It seems likely that Kan knew of this link between homotopy coherence and the comonadic resolutions by at least 1980, (\cf \cite{D&K80a}),  but the construction  does not seem to appear in his work with Dwyer as being linked with coherence until much later. Cordier made the link explicit in \cite{cordier82} and showed how Leitch and Segal's work fit  the pattern.  

Vogt proved in \cite{vogt:1973}, for $\mathcal{C}=Top$, that there was a \emph{homotopy} category of homotopy coherent diagrams of type $\mathsf{A}$ in $Top$\footnote{The definition and composition of homotopy coherent morphisms is difficult to describe briefly.}, and that this was equivalent to the category obtained from $\mathcal{C}^\mathsf{A}$ by inverting the `levelwise' homotopy equivalences, that is those morphisms, $f:X\to Y$, in $\mathcal{C}^\mathsf{A}$ such that for each object, $a$, in $\mathsf{A}$, $f(a)$ is a homotopy equivalence. This result was extended by Cordier and Porter,  \cite{cordierporter:Vogt:1986}, to suitable examples of locally Kan simplicial categories. This gives a number of useful links between the explicit definition of homotopy coherent maps and their interpretation in more homotopical terms.
\subsection{Quasi-categories: once over lightly}\label{qcat}

This  homotopy coherent nerve as defined by  Cordier  is a \emph{quasi-category}, so let us `recall' the definition of that term.

for each object
In a Kan complex, $K$, every horn, $h: \Lambda[n,k]\to K$, with $0\leq k\leq n$, has a filler, that is there is a $n$-simplex whose faces match with those of $h$, where that makes sense, \ie for all faces except the $k^{th}$ one, which is, of course, missing from $\Lambda[n,k]$.  The nerve of any groupoid is a Kan complex, but, in the nerve of a category, \emph{outer} horns, that is those in which $k=0$ or $k=n$, may not have fillers.

 Simplicial sets with the property that all inner horns have fillers were originally called weak Kan complexes, but now they are more often called \emph{quasi-categories}, following the terminology of Joyal in \cite{joyal:2008}.  They are sometimes  called just $\infty$-categories, but that, sometimes, seems lacking in precision, although is a useful term in the context.. 
 
Quasi-categories are one model  of $\infty$-categories. Given any quasi-category, $K$, and, thinking of this object as an $\infty$-category, we define an  \emph{$\infty$-functor} from $K$ to $\mathcal{C}$ to be simply a simplicial map from $K$ to $Ner_{h.c.}(\mathcal{C})$. 
\subsection{Local systems in this context}
We can now return to the ideas of homotopy coherent representations, and generalisations of local systems, so let us take stock, repeating some points from earlier.
Classically a local system on a space, $X$, is a functor from the fundamental groupoid, $\Pi_1(X)$ of $X$ to some category, $\mathcal{A}$, usually this would be the category of sets,  groups, finite groups, Abelian groups, or of vector spaces over a given field.  The idea originated in work on the cohomology of polyhedra, using cochains on their universal covers, see Reidemeister, \cite{Reidemeister:Topologie:1938}, with its full classical development due to Steenrod, \cite{Steenrod:homology-coeffs:1943}, in 1943.  We will use an $\infty$-categorical interpretation and extension of this.

\begin{dfn} (Lurie, \cite{lurie:alg-K-manifold-top:2014}, Lecture 21)  Let $X$ be a topological space and $\mathcal{C}$ be an $\infty$-category. A \emph{local system on $X$ with values in $\mathcal{C}$}\index{local system on a space with values in an $\infty$-category} is a map of simplicial sets,
$$Sing(X)\to \mathcal{C},$$where $Sing(X)$ is the singular simplicial set of $X$. 
The collection of all local systems on $X$ with values in $\mathcal{C}$ gives an $\infty$-category, $Fun(Sing(X),\mathcal{C})$, which is sometimes simply written as $\mathcal{C}^X$.\end{dfn}
It should be said that perhaps this use by Lurie of the term `local system' might better be replaced by `$\infty$-local system', and we will adopt that terminology.  We note that other authors do use other notation such as $\mathsf{Loc}^{\mathcal C}(\pi_{\infty}X)$.  This latter use incorporates the idea that $Sing(X)$ is the `infinity groupoid' of $X$, hence the notation $\pi_{\infty}X$; see for instance \cite{block-smith:higher-Adv}, in which, in the context of a smooth manifold, $M$,  the term `smooth infinity groupoid' is used for the smooth analogue of $Sing(M)$.
\subsection{The link with the classical notion}

\textbf{Example:}  Suppose $C$ is a (small) category, and $\mathcal{C}= Ner(C)$ is the corresponding quasi-category /  $\infty$-category. Let $f:Sing(X)\to Ner(C)$ be a simplicial map, then 
\begin{itemize}
\item for each point $x\in X$, thought of as a singular 0-simplex, $x:\Delta^0\to X$, we have $f(x)$ is in $Ner(C)_0$, so is an object in the category, $C$;
\item for each $a:x\to y$, $a:\Delta^1\to X$, in $Sing(X)_1$, (so $a$ is a path from $x$ to $y$ in $X$), $f(a)$ is a 1-simplex in $Ner(C)$, so is an arrow from $f(x)$ to $f(y)$ in $C$;
\item for each $\sigma\in Sing(X)_2$, we have that the three faces: $f(d_0\sigma)$, $f(d_1\sigma)$, and $f(d_2\sigma)$, form a 2-simplex in $Ner(C)$, so $f(d_0\sigma)\circ f(d_2\sigma)=f(d_1\sigma)$.
\end{itemize}
We will stop there as that will suffice for the moment.  If $a, b: x\to y$ are homotopic paths in $X$ then, from the homotopy, we can build a 2-simplex, $\sigma$ with $d_0(\sigma)=id_y$, $d_1(\sigma)=b$ and $d_2(\sigma)=a$:
$$\xymatrix{&.\ar[dr]^=\ar@{}[d]|>>>>\sigma&\\.\ar[ur]^a\ar[rr]_b&&.}$$
but then $f(a)=f(b)$.  We conclude that $f:Sing(X)\to Ner(C)$ corresponds to a classical local system $\Pi_1(X)\to C$.

\subsection{Variations on the theme of a homotopy coherent nerve}

These constructions can be adapted for dg-categories, as was done by Lurie, \cite{lurie:higher:2017}. They can further be adapted to give a homotopy coherent nerve of an $A_\infty$-category as in Faonte, \cite{faonte:simplicial:TAC:2017}, and we will meet a variant of this in the next section.

We  let $\mathcal{C}$ be a differential graded category and form a simplicial `class', (\ie a possibly `large simplicial set'), $N_{dg}(\mathcal{C})$, which will be called the \emph{differential graded nerve} of $\mathcal{C}$.

The \emph{idea} of $N_{dg}(\mathcal{C})$ is thus as follows:
\begin{itemize}
\item a 0-simplex is simply an object of $\mathcal{C}$;
\item a 1-simplex of $N_{dg}(\mathcal{C})$ is a morphism of $\mathcal{C}$, but that is a pair of objects, $X_0,X_1$, and some element $f\in \mathcal{C}(X_0,X_1)_0$ satisfying  $df=0$;
\item a 2-simplex consists of a triple of objects, $X_0,X_1, X_2$, of $\mathcal{C}$, a triple of maps $f_{0,1}\in \mathcal{C}(X_0,X_1)_0$, $f_{1,2}\in \mathcal{C}(X_1,X_2)_0$ and $f_{0,2}\in \mathcal{C}(X_0,X_2)_0$ satisfying that each $df_{i,j}=0$, and, moreover, an element $f_{0,1,2}\in \mathcal{C}(X_0,X_2)_1$ with $df_{0,1,2} =  (f_{1,2}\circ f_{0,1})-f_{0,2}$, so $f_{0,1,2}$ `fills' the empty 2-simplex made by the three  $f_{i,j}$s ;\\and so on.
\end{itemize}
The exact formulation generalising this description of the bottom layers of this nerve  can be found in Lurie's \cite{lurie:higher:2017}. We record the following.
\begin{thm}\label{dg-nerve is quasicat}
Let $\mathcal{C}$ be a dg-category, then $N_{dg}(\mathcal{C})$ is a quasi-category.\hfill$\square$
\end{thm}
In Faonte, \cite{faonte:simplicial:TAC:2017}, one finds a second way of defining this, which is based on the idea of finding a dg-categorical version of Cordier's $S[n]$. One could take $S[n]$ and convert it into a dg-category by forming the corresponding chain complexes, but using ideas from the Eilenberg-Zilber theorem, we can build a much smaller model for each $n$. In \cite{faonte:simplicial:TAC:2017}, Faonte  specifies this as an $A_\infty$-category, noting that it is in fact a dg-category.  We will briefly sketch his construction. We will work over a commutative ring $\mathbbm{k}$.

The first step is to build $A_\infty$-models for each standard simplex, $\Delta[n]$.  In fact, these are dg-analogues generated by the categories, $[n]= \{0<\ldots < n\}$.  There are no higher compositions. The  $m_k$ for $k\geq 3$ are trivial, so we \emph{do} get actual dg-categories. 

\begin{dfn}The $A_\infty$-category generated by $[n]$, denoted $A_\infty[n]$, is the $A_\infty$-category defined by \begin{itemize}
\item $Ob(A_\infty[n]) = \{0,1, \ldots, n\}$;
\item for $0\leq i,j\leq n$,
$$A_\infty[n](i,j) =\left\{\begin{array}{ll}
\mathbbm{k}\cdot(i,j)& \textrm{if } i\leq j\\
0& \textrm{if } i>j
\end{array}\right.$$
where $\mathbbm{k}\cdot(i,j)$ is the free dg-$\mathbbm{k}$-module on the single element $(i,j)$, and is concentrated in degree 0;
\item the `multiplications', $m_k$, are given by
\begin{itemize}\item $m_1=0$ (by necessity!);
\item $m_2((j,k),(i,j))= (i,k)$;
\item $m_k=0$ if $k\geq 3$.
\end{itemize}
\end{itemize}\end{dfn}

We can now describe the \emph{small} dg-nerve of a dg-category.  We let $\mathsf{A}$ be a dg-category.

\begin{dfn}The \emph{dg-nerve} of $\mathsf{A}$ is the simplicial set, whose $n$-simplices are given by
$$N_{dg}(\mathsf{A})_n= dg\!-\!Cat_k(A_\infty[n],\mathsf{A}),$$
and with face and degeneracy maps induced from the obvious coface and codegeneracy morphisms between the $A_\infty[n]$s.\end{dfn}

If $\mathsf{A}$ is an $A_\infty$-category, as in section \ref{A-infty cat}, almost the same definition works to give an $\mathsf{A}$-nerve.

\begin{dfn}The \emph{$A_\infty$-nerve}\index{A-infinity nerve@$A_\infty$-nerve} of $\mathsf{A}$ is the simplicial set, whose $n$-simplices are given by
$$N_{A_\infty}(\mathsf{A})_n= A_\infty\!-\!Cat_k(A_\infty[n],\mathsf{A}).$$\end{dfn}
We have that Proposition 2.2.12 of Faonte, \cite{faonte:simplicial:TAC:2017}, gives, if
 $\mathsf{A}$ is an $A_\infty$-category, then $N_{A_\infty}(\mathsf{A})$ is a quasi-category.  
 \subsection{Differential graded and $A_\infty$-local systems}
This means one can use Lurie's definition of `local system' with values in an $\infty$-category to obtain compatible meanings for terms such as `dg-local system' or `$A_\infty$-local system'.

This is fine.  It means, for instance, that we can talk about a dg-local system on a space, $X$, as being simply a simplicial mapping from $Sing(X)$ to $N_{dg}(\mathsf{A})$, for whatever dg-category, $\mathsf{A}$, that we need to use.  That, however, ignores the potential of the basic idea.  We can replace $Sing(X)$ by other simplicial sets as necessary.  This idea has been exploited by various people as we will outline next. We will examine some of them more closely. All these are instances of what may be called representations up to coherent homotopy. All of them correspond to some form of fibration or fibred categorical interpretation. \label{Block-Smith}

If we replace $Sing(X)$ by an arbitrary simplicial set, $K$,   then  we get the notion of dg-local systems and $A_\infty$-local systems introduced by Block and Smith, \cite{block-smith:higher-Adv}.  There they comment that these $A_\infty$-local systems  are closely related to the notion of $A_\infty$-functors introduced by Igusa, \cite{igusa:twisting}.  That latter paper links back to twisted tensor products and in  \cite{RZ:ttp}, Rivera and Zeinlein link colimits of a dg-local  system of chain complexes to those same twisted tensor products.  Here it is important to note that the `colimit' is to be taken in the homotopy  coherent sense. Also it is worth noting that $N_{dg}$ has a left adjoint, again interpreted at the appropriate level of  $\infty$-categorification.

\section{Good open covers, fibrations and infinity local systems}
To tie this slightly abstract theory into more elementary topics, we return to the classical subject of topological fibrations and assume given a fibration, $p:E\to B$, where $B$ has a good open cover, $\mathcal{U}$.  Recall that an open cover, $\mathcal{U}=\{U_\alpha\}$, is said to be `good' if  all the $U_\alpha$ as well as all their finite intersections are contractible or empty. We note that every paracompact smooth manifold admits a good open cover.  In fact, it admits a differentiable good open cover, (see \cite{FSS:Cech:2012}), \ie one in which all non-empty finite  intersections are diffeomorphic to an open ball, although we will not use smoothness here, merely existence. Similarly any triangulation of a manifold gives a good open cover, and the nerve of that cover will be homeomorphic to the given manifold.

If $\mathcal{U}$ is a good cover, then, for each family, $\sigma=\{\alpha_0,\ldots, \alpha_n\}$, of indices, and  writing $U_\sigma=\bigcap_{i=0}^nU_{\alpha_i}$, we choose an (explicit) homotopy equivalence between $U_\sigma$ and a point, $x_\sigma\in U_\sigma$. We write $F_\sigma$ for $p^{-1}(x_\sigma)$, and obtain an explicit  homotopy equivalence between $p^{-1}(U_\sigma)$ and $F_\sigma$.

We now need a variant of a result that has, in various shapes and forms, been used, behind the scenes, at several points earlier in this survey.  As we mentioned on page \pageref{replacement}, Vogt proved in \cite{vogt:1973} that given a small category, $\mathsf{A}$, there was an equivalence of categories between $Ho(Top^\mathsf{A})$, obtained from $Top^\mathsf{A}$ by inverting levelwise homotopy equivalences, and the homotopy category of homotopy coherent diagrams $Coh(\mathsf{A},Top)$. In \cite{cordierporter:Vogt:1986}, Cordier and TP extended this to considering  $\mathcal{C}^\mathsf{A}$, where $\mathcal{C}$ is a locally Kan simplicially enriched category (with some extra conditions).  One of the key results in those theories is to show that given a diagram, $X:\mathsf{A}\to \mathcal{C}$, then specifying some homotopy equivalences, 
$$\xymatrix{X(a)\ar[r]<.5ex>^{h_a}&Y(a)\ar[l]<.5ex>^{k_a}},$$ together with the homotopies $h_ak_a\simeq 1_{Y(a)}$ and $k_ah_a\simeq 1_{X(a)}$, for each object $a$ in $\mathsf{A}$, it is easy to build a homotopy coherent diagram based on the $Y(a)$s; see \cite{C&P88} for more detail.  That result, or rather variants or special cases  of it in the dg- or $A_\infty$-setting, had been used many times in the literature, and was the idea behind some of Segal's constructions in \cite{segal:1974}.  We will call it the `replacement result' for convenience.\label{replacement}

Going back to the good open cover, $\mathcal{U}$, we take $Face(B, \mathcal{U})$ to be the face poset of the \v{C}ech nerve of $\mathcal{U}$, so $Face(B, \mathcal{U})$ has the $\sigma$s for which $U_\sigma$ is non-empty as its elements with $\leq$ being $\subseteq$.  The nerve of this poset is the barycentric subdivision of the  \v{C}ech nerve of $\mathcal{U}$. Thinking of this poset as a category, we have a functor,
$\mathcal{F}: Face(B, \mathcal{U})\to Top$, sending $\sigma$ to $p^{-1}(U_\sigma)$, and, using the replacement result, we get a homotopy coherent diagram,$$\mathcal{F}:K\to Ner_{h.c.}(Top),$$ where $K$ is $Ner(Face(B, \mathcal{U}))$.  This is an $\infty$-local system derived from $p:E\to B$ and is determined by it up to homotopy.  This is one way to extend to the infinity categorical setting, the classical action of the fundamental groupoid of the base on the sets of path components of the fibres. We note that $K$ is a quasi-category, but not usually an infinity groupoid. The homotopy transition cocycles of Wirth's thesis, \cite{wirth:diss}, and \cite{jw-jds}, are another manifestation of this type of idea.

If one applies many of the usual functors from $Top$ to useful algebraic categories, for instance, of partial algebraic models for homotopy types, it is easy to convert this $\infty$-local system into other forms \eg dg-, $A_\infty$- or groupoid valued ones.  This is closely related to the non-abelian cohomological setting that we mentioned in section \ref{non-abelian}, but we will not follow up on that connection here.

We will also not explore reconstruction techniques via gluings or analogues of the Grothendieck construction / homotopy colimit, but note that, as we are almost explicitly using the Joyal -  Lurie theory of quasi-categories here, we have at our disposal the tools of that theory, including Lurie's straightening /  unstraightening methodology, \cf \cite{lurie:HighTop2006}, which is a form of $\infty$-categorical Grothendieck construction, to investigate the analogues of reconstruction in this setting.

\section{Postscript}

There are various topics that deserved to be mentioned, but due to limitations of space and time, have had to be omitted. Some of these were already mentioned in our earlier survey article, \cite{TP-JDS-hcr:2022}.  We would draw particular notice to possible higher dimensional categorified versions of flat connections, holonomy, parall el transport  and the Riemann-Hilbert correspondence. These were almost touched on in section \ref{Block-Smith}, where the use by Block and Smith of $A_\infty$-local systems, in \cite{block-smith:higher-Adv},  was noted, but these would ideally have been discussed more fully. Another topic would have been to link Lurie's local systems with the actions of path spaces and the coherence problems that arise there. Links with non-abelian cohomology might also have been expanded, but all these topics will have to wait for another time. There are also more speculative areas including extensions of the ideas to stratified spaces, orbifolds and related structures, but that would have required a much longer survey. We should also point out that homotopy coherence has  appeared in the math-physics literature, especially in the guise of $L_\infty$-structures, for instance in gauge field theory, topics that we have not touched on here.

This survey is meant as an appetiser, so we have tried not to have a bibliography longer that the body of the paper.  By necessity this meant that some people may feel aggrieved not to be mentioned, and some topics have not even been treated at all, but `such is life'.

\section*{Acknowledgements.} The contents of this survey article have been much influenced by the group of researchers, who have discussed the area of infinity local systems, homotopy coherence, etc., over Zoom in the last six months. They include Manuel Rivera, Mahmoud Zeinalian, Francis Bischoff, Jo\~{a}o Faria Martins, Domenico Fiorenza, Camilo Arias Abad and other more occasional visitors.

\bibliographystyle{elsarticle-harv} 

\bibliography{Fields}
\ddo